# Response Matrix Benchmark for the 1D Transport Equation with Matrix Scaling

B.D. Ganapol[1]
University of Arizona
Department of Aerospace and Mechanical Engineering
Tucson AZ, USA

J. K. Patel
University of Michigan
Nuclear Engineering and Radiological Sciences Department
Ann Arbor, MI, USA

**ABSTRACT**

The linear 1D transport equation is likely the most solved transport equation in radiative transfer and neutron transport investigations. Nearly every method imaginable has been applied to establish solutions, including Laplace and Fourier transforms, singular eigenfunctions, solutions of singular integral equation, PN expansions, double PN expansions, Chebychev expansions, Lagrange polynomial expansions, numerical discrete ordinates with finite difference, analytical discrete ordinates, finite elements, solutions to integral equations, adding and doubling, invariant imbedding, solution of Ricatti equations and response matrix methods-- and probably more methods of which the authors are unaware. Of those listed, the response matrix solution to the discrete ordinates form of the 1D transport equation is arguably the simplest and most straightforward. Here, we propose another response of exponential solutions but to the first order equation enabled by matrix scaling.

## 1. INTRODUCTION

The neutron transport and radiative transfer equations offer a rich variety of opportunities for mathematical and numerical innovation. The opportunities range from discrete ordinate approximations, to singular eigenfunction expansions, even in 3D. With time, new solutions seem to continuously appear in the literature. In this presentation, we re-visit the Response Matrix/Discrete Ordinate Method (RM/DOM) proposed several years ago [1] by the first author, which was based on the second order even/odd parity form of the transport equation coupled to Wynn-epsilon acceleration. The RM/DOM was successful in generating a 7-place solution for highly anisotropic scattering such as the 300-term CloudC1 kernel. The emphasis here will again be highly forward peaked scattering, however, using the first order transport equation itself. One enables the linear algebraic solution to the discrete ordinates equations, represented as an exponential in first order form by including the partitioned identity matrix as an intermediate step. The result is the stabilization of the solution for large eigenvalues. An example comparison of the first order RM/DOM (called 1st/RM/DOM) with RM/DOM solved in second order form (2nd/RM/DOM) confirms all 7 places and nearly 8-places for the CloudC1 scattering phase function.

In today's transport theory applications, 1D methods are hardly at the forefront when compared to 3D unstructured mesh solutions one finds in the literature and at conferences. However, the sophisticated

---

reader is aware of the role 1D benchmarks play in transport theory. In particular, 1D has many of the features of 2D or 3D transport, but is more readily available by emphasizing scattering and streaming rather than geometric complexity. Combine simplicity with high precision and the notion if an advanced transport algorithm cannot achieve a desired precision in 1D, then regardless of how advanced the method, it must be viewed with caution. Thus, the value of a 1D benchmark and the motivation for pursuing yet another in this presentation.

## 2. THE FIRST ORDER SOLUTION

### 2.1. The Angularly Discretized 1D Transport Equation

In a 1D homogeneous slab, the transport equation for the intensity $I(\tau,\mu)$ (or flux) satisfies

$$\left[\mu\frac{\partial}{\partial\tau}+1\right]I(\tau,\mu)=\omega\int_{-1}^{1}d\mu' f(\mu',\mu)I(\tau,\mu'),\ \tau_0<\tau<\tau_1 \quad (1a)$$

for particle motion in direction (cosine) $\mu$ at position (optical depth) $\tau$. The single scatter albedo is $\omega$ and $f(\mu',\mu)$ is the scattering phase function. The incoming flux on the free surfaces is

$$\begin{aligned}I(\tau_0,\mu)&=f(\mu)\\ I(\tau_1,-\mu)&=g(\mu)\end{aligned} \quad (1b,c)$$

for $0\leq\mu\leq 1$. A Gauss quadrature, for $N$ even, gives the following discretization of directions:

$$\mu=\begin{cases}-\mu_m\\ \mu_{N+m}\equiv\mu_m\end{cases},\ m=1,...,N,\ \mu_m>0, \quad (2a)$$

to approximate the collision integral, and Eq(1a) becomes in a given direction $\mu_m$

$$\left[\mu_m\frac{\partial}{\partial\tau}+1\right]I(\tau,\mu_m)=\omega\sum_{m'=1}^{2N}\alpha_{m'}f(\mu_{m'},\mu_m)I(\tau,\mu_{m'}). \quad (2b)$$

Now consider particle motion in directions toward the far surface (+) at $x = \tau_1$ and toward the near surface (−) at $x = \tau_0$ as shown in Fig. 1. In vector form, Eq(2b) becomes [1]

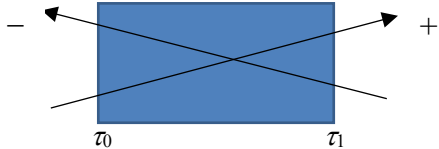

Fig. 1. Directions of particle motion.

$$\frac{d\mathbf{I}^{\pm}(\tau)}{d\tau}=\mp\mathbf{M}^{-1}\left(\mathbf{I}_N-\mathbf{C}^{\pm\pm}\mathbf{W}\right)\mathbf{I}^{\pm}(\tau)\pm\mathbf{M}^{-1}\mathbf{C}^{\pm\mp}\mathbf{W}\mathbf{I}^{\mp}(\tau), \quad (3a\pm)$$

and in a more compact vector form

$$\frac{d\boldsymbol{I}^+(\tau)}{d\tau} = \boldsymbol{\beta}\boldsymbol{I}^-(\tau) - \boldsymbol{\alpha}\boldsymbol{I}^+(\tau) \qquad (3b\pm)$$

$$\frac{d\boldsymbol{I}^-(\tau)}{d\tau} = \boldsymbol{\alpha}\boldsymbol{I}^-(\tau) - \boldsymbol{\beta}\boldsymbol{I}^+(\tau),$$

where

$$\boldsymbol{\alpha} \equiv \boldsymbol{M}^{-1}\left(\boldsymbol{I}_N - \boldsymbol{C}^{++}\boldsymbol{W}\right) \qquad (3c,d)$$

$$\boldsymbol{\beta} \equiv \boldsymbol{M}^{-1}\boldsymbol{C}^{+-}\boldsymbol{W},$$

and $\boldsymbol{I}_N$ is the identity matrix of size $N$. Continuing,

$$\boldsymbol{I}^{\mp}(\tau) \equiv \left[ I_{\genfrac{\{}{\}}{0pt}{}{1}{N+1}}(\tau) \quad I_{\genfrac{\{}{\}}{0pt}{}{2}{N+2}}(\tau) \quad \cdots \quad I_{\genfrac{\{}{\}}{0pt}{}{N}{2N}}(\tau) \right]^T \qquad (4a\mp)$$

$$\boldsymbol{M} \equiv diag\{\mu_m\}$$
$$\boldsymbol{W} \equiv diag\{\alpha_m\} \qquad (4b,c)$$

$$\boldsymbol{C}^{\pm\pm}\boldsymbol{W} \equiv \{\omega\alpha_m f(\pm\mu_{m'}, \pm\mu_m); m' = 1,...,N; m = 1,...,N\}$$
$$\boldsymbol{C}^{\pm\mp}\boldsymbol{W} \equiv \{\omega\alpha_m f(\pm\mu_{m'}, \mp\mu_m); m' = 1,...,N; m = 1,...,N\}. \qquad (4d^{\pm\pm},e^{\pm\mp})$$

From scattering symmetry

$$f(\mu_i, -\mu_j) = f(-\mu_i, \mu_j)$$
$$f(-\mu_i, -\mu_j) = f(\mu_i, \mu_j). \qquad (5)$$

Note, boldface indicates either a vector or matrix distinguished by context. Finally, the equation to solve from Eqs(3a±) is the first order ODE

$$\frac{d\boldsymbol{I}(\tau)}{d\tau} = \boldsymbol{A}\boldsymbol{I}(\tau) \qquad (6a)$$

with

$$\boldsymbol{I}(\tau) \equiv \begin{bmatrix} \boldsymbol{I}^+(\tau) \\ \boldsymbol{I}^-(\tau) \end{bmatrix} \qquad (6b)$$

$$\boldsymbol{A} \equiv \begin{bmatrix} -\boldsymbol{\alpha} & \boldsymbol{\beta} \\ -\boldsymbol{\beta} & \boldsymbol{\alpha} \end{bmatrix}. \qquad (6c)$$

### 2.2. First order solution

For a constant $\boldsymbol{A}$ square matrix of size $2N$, the formal solution to Eqs(6) is

$$I(\tau) \equiv e^{A\tau} I(\tau_0). \tag{7}$$

Therefore, the spatial representation is analytical and requires the evaluation of the matrix exponential, which is numerically considered a formidable task [2] with many variations. For *N*, say under 500, a suitable evaluation, via matrix diagonalization, follows if the matrix has independent eigenvectors as is true here. Diagonalization is more direct than the alternative of specifying eigenvectors and assuming a complete set to form a spectral expansion. In any case, as can be shown, spectral expansion is a consequence of diagonalization.

Diagonalization of *A* gives

$$\boldsymbol{A} = \boldsymbol{T}\boldsymbol{\lambda}\boldsymbol{T}^{-1}, \tag{8a}$$

where the eigenvalues are real symmetric sets of negative and positive values with negative/positive ordering

$$\begin{aligned}\tilde{\lambda}^- &\equiv \left[-\lambda_j; j = 1,\ldots,N\right] \\ \tilde{\lambda}^+ &\equiv \left[\lambda_j; j = N+1,\ldots,2N\right].\end{aligned} \tag{8b}$$

Ordering of the eigenvalues is generally arbitrary; but here, we group the negatives and positives together in ascending order, which seems reasonable. Then, combining the sets into a diagonal matrix, the eigenvalue matrix is

$$\lambda \equiv \begin{bmatrix} \tilde{\lambda}^- & \mathbf{0} \\ \mathbf{0} & \tilde{\lambda}^+ \end{bmatrix}. \tag{9}$$

*T* is the eigenvector matrix, where each column corresponds to an eigenvalue of the diagonal eigenvalue matrix Eq(9).

The eigenvectors, by column, inherit the ordering of the eigenvalues,

$$\boldsymbol{T} \equiv \begin{bmatrix} \boldsymbol{T}_1 & \boldsymbol{T}_2 & \ldots & \boldsymbol{T}_{2N} \end{bmatrix}^T = \begin{bmatrix} \boldsymbol{T}^- & \boldsymbol{T}^+ \end{bmatrix}^T \tag{10a}$$

with

$$\begin{aligned}\boldsymbol{T}^- &\equiv \begin{bmatrix} \boldsymbol{T}_1 & \boldsymbol{T}_2 & \ldots & \boldsymbol{T}_N \end{bmatrix}^T = \begin{bmatrix} \boldsymbol{T}_N^- & \boldsymbol{T}_{N-1}^- & \ldots & \boldsymbol{T}_1^- \end{bmatrix}^T \\ \boldsymbol{T}^+ &\equiv \begin{bmatrix} \boldsymbol{T}_{N+1} & \boldsymbol{T}_{N+2} & \ldots & \boldsymbol{T}_{2N} \end{bmatrix}^T = \begin{bmatrix} \boldsymbol{T}_1^+ & \boldsymbol{T}_2^+ & \ldots & \boldsymbol{T}_N^+ \end{bmatrix}^T\end{aligned} \tag{10b$\mp$}$$

and from Eqs(9) and (10), *A* becomes

$$A = T \begin{bmatrix} \tilde{\lambda}^- & 0 \\ 0 & \tilde{\lambda}^+ \end{bmatrix} T^{-1}. \tag{10c}$$

The diagonalized matrix exponential is therefore

$$e^{A\tau} = T \begin{bmatrix} \Gamma^-(\tau) & 0 \\ 0 & \Gamma^+(\tau) \end{bmatrix} T^{-1}, \tag{11a}$$

with

$$\Gamma^-(\tau) \equiv diag\{e^{-\lambda_{N-k+1}\tau}; k=1,...,N\}$$
$$\Gamma^+(\tau) \equiv diag\{e^{\lambda_k \tau}; k=1,...,N\}. \tag{11∓}$$

Therefore, the algebraic solution is simply

$$I(\tau) = T \begin{bmatrix} \Gamma^-(\tau) & 0 \\ 0 & \Gamma^+(\tau) \end{bmatrix} T^{-1} I(\tau_0). \tag{12}$$

If we partition $T$ and $T^{-1}$ into four square matrices each of size $N$

$$T \equiv \begin{bmatrix} T_1 & T_2 \\ T_3 & T_4 \end{bmatrix} \tag{13a}$$

$$T^{-1} \equiv \begin{bmatrix} T_1 & T_2 \\ T_3 & T_4 \end{bmatrix}^{-1}, \tag{13b}$$

the matrix exponential becomes

$$e^{A\tau} = \begin{bmatrix} T_1 & T_2 \\ T_3 & T_4 \end{bmatrix} \begin{bmatrix} \Gamma^-(\tau) & 0 \\ 0 & \Gamma^+(\tau) \end{bmatrix} \begin{bmatrix} T_1 & T_2 \\ T_3 & T_4 \end{bmatrix}^{-1}. \tag{14}$$

Finally, Eq(12) is

$$I(\tau) = \begin{bmatrix} T_1 & T_2 \\ T_3 & T_4 \end{bmatrix} \begin{bmatrix} \Gamma^-(\tau)_1 & 0 \\ 0 & \Gamma^+(\tau)_4 \end{bmatrix} \begin{bmatrix} T_1 & T_2 \\ T_3 & T_4 \end{bmatrix}^{-1} I(\tau_0). \tag{15}$$

Note the exponential diagonals are further identified by their quadrant, either 1 or 4. Therein lies the difficulty of evaluating Eq(15). Since some eigenvalues can be large for $\omega$ positive near unity, at least one element of $\Gamma^+(\tau)_4$ will give an overflow. To eliminate this possibility, we introduce scaling. This issue has be been noted previously in Siewert's ADO method [3].

## 3. SCALING

### 3.1 Including the identity matrix
To avoid overflow, we insert a scaling factor through the identity $I_{2N}$ (size 2N) factored as

$$\begin{bmatrix} \Gamma^-(\tau)_1 & 0 \\ 0 & \Gamma^+(\tau)_4 \end{bmatrix} = \begin{bmatrix} \Gamma^-(\tau)_1 & 0 \\ 0 & \Gamma^+(\tau)_4 \end{bmatrix} \underbrace{\begin{bmatrix} I & 0 \\ 0 & \Gamma^-(\tau_1)_4 \end{bmatrix} \begin{bmatrix} I & 0 \\ 0 & \Gamma^-(\tau_1)_4 \end{bmatrix}^{-1}}_{I_{2N}(\text{size }2N)}, \quad (16a)$$

to give

$$\begin{bmatrix} \Gamma^-(\tau)_1 & 0 \\ 0 & \Gamma^+(\tau)_4 \end{bmatrix} = \begin{bmatrix} \Gamma^-(\tau)_1 & 0 \\ 0 & \Gamma^-(\tau_1-\tau)_4 \end{bmatrix} \begin{bmatrix} I & 0 \\ 0 & \Gamma^-(\tau_1)_4 \end{bmatrix}^{-1}, \quad (16b)$$

since

$$\Gamma^+(\tau)_4 \Gamma^-(\tau_1)_4 = \Gamma^-(-\tau)_4 \Gamma^-(\tau_1)_4 = \Gamma^-(\tau_1-\tau)_4. \quad (16c)$$

Thus, when Eq(16b) is introduced into Eq(15),

$$I(\tau) = \begin{bmatrix} T_1 & T_2 \\ T_3 & T_4 \end{bmatrix} \begin{bmatrix} \Gamma^-(\tau)_1 & 0 \\ 0 & \Gamma^-(\tau_1-\tau)_4 \end{bmatrix} \begin{bmatrix} I & 0 \\ 0 & \Gamma^-(\tau_1)_4 \end{bmatrix}^{-1} \begin{bmatrix} T_1 & T_2 \\ T_3 & T_4 \end{bmatrix}^{-1} I(\tau_0), \quad (17a)$$

with inverses combined and matrices multiplied, Eq(15) becomes

$$I(\tau) = \begin{bmatrix} T_1\Gamma^-(\tau)_1 & T_2\Gamma^-(\tau_1-\tau)_4 \\ T_3\Gamma^-(\tau)_1 & T_4\Gamma^-(\tau_1-\tau)_4 \end{bmatrix} \begin{bmatrix} T_1 & T_2\Gamma^-(\tau_1)_4 \\ T_3 & T_4\Gamma^-(\tau_1)_4 \end{bmatrix}^{-1} I(\tau_0). \quad (17b)$$

resulting in the elimination of all positive exponents. If we define the unknown vector

$$\begin{bmatrix} \alpha \\ \beta \end{bmatrix} \equiv \begin{bmatrix} T_1 & T_2\Gamma^-(\tau_1)_4 \\ T_3 & T_4\Gamma^-(\tau_1)_4 \end{bmatrix}^{-1} \begin{bmatrix} I^+(\tau_0) \\ I^-(\tau_0) \end{bmatrix}, \quad (18a)$$

then, Eq(17b) is

$$I(\tau) = \begin{bmatrix} T_1\Gamma^-(\tau)_1 & T_2\Gamma^-(\tau_1-\tau)_4 \\ T_3\Gamma^-(\tau)_1 & T_4\Gamma^-(\tau_1-\tau)_4 \end{bmatrix}\begin{bmatrix}\alpha\\\beta\end{bmatrix} \quad (18b)$$

or by component is

$$I^+(\tau) = T_1\Gamma^-(\tau)_1\alpha + T_2\Gamma^-(\tau_1-\tau)_4\beta$$
$$I^-(\tau) = T_3\Gamma^-(\tau)_1\alpha + T_4\Gamma^-(\tau_1-\tau)_4\beta. \quad (18c)$$

By introducing $\tau = \tau_1$, and $\tau_0$, Eqs(18c) become respectively,

$$\begin{bmatrix}I^+(\tau_1)\\I^-(\tau_0)\end{bmatrix} = \begin{bmatrix} T_1\Gamma^-(\tau_1)_1 & T_2 \\ T_3\Gamma^-(\tau_0)_1 & T_4\Gamma^-(\tau_1-\tau_0)_4 \end{bmatrix}\begin{bmatrix}\alpha\\\beta\end{bmatrix}. \quad (19)$$

Next, we introduce $\tau = \tau_0$, and $\tau_1$ into Eqs(18c) respectively,

$$I^+(\tau_0) = T_1\Gamma^-(\tau_0)_1\alpha + T_2\Gamma^-(\tau_1-\tau_0)_4\beta$$
$$I^-(\tau_1) = T_3\Gamma^-(\tau_1)_1\alpha + T_4\beta \quad (20a)$$

and solve for the coefficient vectors $\alpha, \beta$

$$\begin{bmatrix}\alpha\\\beta\end{bmatrix} = \begin{bmatrix} T_1\Gamma^-(\tau_0)_1 & T_2\Gamma^-(\tau_1-\tau_0)_4 \\ T_3\Gamma^-(\tau_1)_1 & T_4 \end{bmatrix}^{-1}\begin{bmatrix}I^+(\tau_0)\\I^-(\tau_1)\end{bmatrix}. \quad (20b)$$

When Eq(20b) is introduced into Eqs(19), we find

$$\begin{bmatrix}I^+(\tau_1)\\I^-(\tau_0)\end{bmatrix} = R\begin{bmatrix}I^+(\tau_0)\\I^-(\tau_1)\end{bmatrix}, \quad (21a)$$

where the response matrix without the troublesome unbounded exponentials emerges as

$$R = \begin{bmatrix} T_1\Gamma^-(\tau_1)_1 & T_2 \\ T_3\Gamma^-(\tau_0)_1 & T_4\Gamma^-(\tau_1-\tau_0)_4 \end{bmatrix}\begin{bmatrix} T_1\Gamma^-(\tau_0)_1 & T_2\Gamma^-(\tau_1-\tau_0)_4 \\ T_3\Gamma^-(\tau_1)_1 & T_4 \end{bmatrix}^{-1}. \quad (21b)$$

Therefore, to obtain the outgoing intensities at the slab surfaces, one only need to multiply the response matrix by the incoming vector intensities.

The interior intensities come from Eqs(18b) and (20b)

$$\begin{bmatrix} I^+(\tau) \\ I^-(\tau) \end{bmatrix} = \begin{bmatrix} T_1\Gamma^-(\tau)_1 & T_2\Gamma^-(\tau_1-\tau)_4 \\ T_3\Gamma^-(\tau)_1 & T_4\Gamma^-(\tau_1-\tau)_4 \end{bmatrix} \begin{bmatrix} T_1\Gamma^-(\tau_0)_1 & T_2\Gamma^-(\tau_1-\tau_0)_4 \\ T_3\Gamma^-(\tau_1)_1 & T_4 \end{bmatrix}^{-1} \begin{bmatrix} I^+(\tau_0) \\ I^-(\tau_1) \end{bmatrix}.$$
(22)

The solution of Eq(22), is arguably the most straightforward solution of the discrete ordinates equations based only on the knowledge of the eigenvalues and eigenvectors of matrix $A$. One can easily extend the analysis to heterogeneous media.

### 3.2 A demonstration
For comparison of 1st/RM/DOM with 2nd/RM/DOM intensities, we assume the 300-term CloudC1 phase function [4] for a slab of 64 *mfp* thickness. A perpendicular beam ($\mu_0 = 1$) enters the near surface $\mu > 0$

$$I(0,\mu) = \frac{1}{2}\delta(\mu-1)$$
$$I(\tau_1,-\mu) = 0$$
(24)

with vacuum bordering the far surface. For the single scatter albedo of unity, Table Ia gives a seven-place benchmark by 2nd/DOM/RM for exiting and four asymmetric interior intensities in eleven-edit directions. We make the comparison with 1st/DOM/RM through faux (false) quadratures, where the eleven-edit directions are included in the quadrature list with zero weight. In this way, the transport algorithm itself interpolates the edits but does not affect the solution. To enable convergence of the solution in quadrature order $N$, the quadrature order is incremented until the relative error between edits fall below a prescribed relative error. Both calculations apply a Radau quadrature. All entries between the two calculations agree to all digits quoted for $N_{1st} = 300$, $N_{2nd} = 295$ for 1st and 2nd benchmarks respectively.

**Table Ia**. 7-place benchmark $\tau_1=64$, $\omega = 1$ for both 1st and 2nd RM/DOM for asymmetric spatial edits

| $\mu\backslash\tau$ | 0 | $\tau_1/20$ | $\tau_1/5$ | $\tau_1/2$ | $3\tau_1/4$ | $\tau_1$ |
|---|---|---|---|---|---|---|
| -1.000E+00 | 1.0636984E+00 | 1.0062387E+00 | 8.5824229E-01 | 5.2453336E-01 | 2.4600228E-01 | 0.0000000E+00 |
| -8.000E-01 | 9.5407647E-01 | 9.9828274E-01 | 8.8052819E-01 | 5.4744427E-01 | 2.6883574E-01 | 0.0000000E+00 |
| -6.000E-01 | 8.2471232E-01 | 9.7909867E-01 | 9.0255626E-01 | 5.7035518E-01 | 2.9173427E-01 | 0.0000000E+00 |
| -4.000E-01 | 7.1143850E-01 | 9.5800446E-01 | 9.2437584E-01 | 5.9326601E-01 | 3.1464747E-01 | 0.0000000E+00 |
| -2.000E-01 | 5.5848173E-01 | 9.2583435E-01 | 9.4589447E-01 | 6.1617674E-01 | 3.3756297E-01 | 0.0000000E+00 |
| 0.000E+00 | 2.5158245E-01 | 8.7951999E-01 | 9.6701431E-01 | 6.3908730E-01 | 3.6047812E-01 | 3.9263859E-02 |
| 2.000E-01 | 0.0000000E+00 | 8.1871840E-01 | 9.8770148E-01 | 6.6199770E-01 | 3.8339248E-01 | 8.8948974E-02 |
| 4.000E-01 | 0.0000000E+00 | 7.5459856E-01 | 1.0082266E+00 | 6.8490802E-01 | 4.0630611E-01 | 1.1838392E-01 |
| 6.000E-01 | 0.0000000E+00 | 7.3765495E-01 | 1.0300017E+00 | 7.0781884E-01 | 4.2921921E-01 | 1.4513931E-01 |
| 8.000E-01 | 0.0000000E+00 | 8.8704553E-01 | 1.0589106E+00 | 7.3073267E-01 | 4.5213189E-01 | 1.7054334E-01 |
| 1.000E+00 | 0.0000000E+00 | 8.0745964E+01 | 1.2605960E+00 | 7.5366350E-01 | 4.7504419E-01 | 1.9523120E-01 |

Table Ib attempts an 8-place comparison, which is successful except for two entries (emboldened) for $N_{1st}$ = 350. This comparison also indicates that comparisons in the last digits can give erroneous results. While Table Ib shows two discrepancies, by limiting $N_{1st}$ to 300, we find perfect 8-place agreement.

The method of comparing benchmarks is somewhat arbitrary in that it is just for a snapshot of the independent variables. For example, another choice of spatial positions could change the agreement in the last digit, or not. Table IIa gives the 8-place benchmark (2nd) for uniformly distributed spatial edits. Now we observe perfect agreement for $N_{1st} = 350$.

**Table Ib**. 8-place benchmark (2$^{nd}$) for $\tau_1$=64, $\omega$ = 1 for both 1$^{st}$ and 2$^{nd}$RM/DOM for asymmetric spatial edits

| μ\τ | 0 | $\tau_1$/20 | $\tau_1$/5 | $\tau_1$/2 | 3$\tau_1$/4 | $\tau_1$ |
|---|---|---|---|---|---|---|
| -1.000E+00 | 1.06369836E+00 | 1.00623866E+00 | 8.58242288E-01 | 5.24533358E-01 | 2.46002278E-01 | 0.00000000E+00 |
| -8.000E-01 | 9.54076470E-01 | 9.98282743E-01 | 8.80528187E-01 | 5.47444267E-01 | 2.68835741E-01 | 0.00000000E+00 |
| -6.000E-01 | 8.24712324E-01 | 9.79098672E-01 | 9.02556262E-01 | 5.70355177E-01 | 2.91734272E-01 | 0.00000000E+00 |
| -4.000E-01 | 7.11438505E-01 | 9.58004461E-01 | 9.24375836E-01 | 5.93266014E-01 | 3.14647466E-01 | 0.00000000E+00 |
| -2.000E-01 | 5.58481727E-01 | 9.25834350E-01 | 9.45894466E-01 | 6.16176736E-01 | 3.37562967E-01 | 0.00000000E+00 |
| 0.000E+00 | 2.51582455E-01 | 8.79519989E-01 | 9.67014307E-01 | 6.39087303E-01 | 3.60478121E-01 | 3.92638589E-02 |
| 2.000E-01 | 0.00000000E+00 | 8.18718402E-01 | 9.87701480E-01 | 6.61997697E-01 | 3.83392476E-01 | 8.89489741E-02 |
| 4.000E-01 | 0.00000000E+00 | 7.54598563E-01 | 1.00822662E+00 | 6.84908019E-01 | 4.06306114E-01 | 1.18383923E-01 |
| 6.000E-01 | 0.00000000E+00 | 7.37654954E-01 | 1.03000170E+00 | 7.07818836E-01 | 4.29219210E-01 | 1.45139314E-01 |
| 8.000E-01 | 0.00000000E+00 | 8.87045535E-01 | 1.05891057E+00 | 7.30732673E-01 | 4.52131890E-01 | 1.70543344E-01 |
| 1.000E+00 | 0.00000000E+00 | 8.07459640E+01 | 1.26059598E+00 | 7.53663496E-01 | 4.75044193E-01 | 1.95231197E-01 |

**Table IIa**. 8-place benchmark for $\tau_1$=64, $\omega$ = 1 for both 1$^{st}$ and 2$^{nd}$RM/DOM for uniform spatial edits

| μ\τ | 0 | $\tau_1$/5 | 2$\tau_1$/5 | 3$\tau_1$/5 | 4$\tau_1$/5 | $\tau_1$ |
|---|---|---|---|---|---|---|
| -1.000E+00 | 1.06369836E+00 | 8.58242289E-01 | 6.35970510E-01 | 4.13094529E-01 | 1.90467611E-01 | 0.00000000E+00 |
| -8.000E-01 | 9.54076470E-01 | 8.80528187E-01 | 6.58878508E-01 | 4.36004258E-01 | 2.13094527E-01 | 0.00000000E+00 |
| -6.000E-01 | 8.24712324E-01 | 9.02556262E-01 | 6.81785223E-01 | 4.58915261E-01 | 2.35953990E-01 | 0.00000000E+00 |
| -4.000E-01 | 7.11438505E-01 | 9.24375836E-01 | 7.04690673E-01 | 4.81826511E-01 | 2.58872714E-01 | 0.00000000E+00 |
| -2.000E-01 | 5.58481727E-01 | 9.45894466E-01 | 7.27594473E-01 | 5.04737793E-01 | 2.81800546E-01 | 0.00000000E+00 |
| 0.000E+00 | 2.51582455E-01 | 9.67014307E-01 | 7.50496082E-01 | 5.27649059E-01 | 3.04726966E-01 | 3.92638589E-02 |
| 2.000E-01 | 0.00000000E+00 | 9.87701480E-01 | 7.73395273E-01 | 5.50560298E-01 | 3.27650253E-01 | 8.89489741E-02 |
| 4.000E-01 | 0.00000000E+00 | 1.00822662E+00 | 7.96293473E-01 | 5.73471520E-01 | 3.50570747E-01 | 1.18383923E-01 |
| 6.000E-01 | 0.00000000E+00 | 1.03000170E+00 | 8.19198648E-01 | 5.96382767E-01 | 3.73489122E-01 | 1.45139314E-01 |
| 8.000E-01 | 0.00000000E+00 | 1.05891057E+00 | 8.42146047E-01 | 6.19294219E-01 | 3.96405868E-01 | 1.70543344E-01 |
| 1.000E+00 | 0.00000000E+00 | 1.26059598E+00 | 8.65364223E-01 | 6.42206808E-01 | 4.19321107E-01 | 1.95231197E-01 |

**Table IIb**. 8-place benchmark for $\tau_1$=64, $\omega$ = 1 for both 1$^{st}$ and 2$^{nd}$RM/DOM for uniform spatial edits and 21 directional edits

| μ\τ | 0 | $\tau_1$/5 | 2$\tau_1$/5 | 3$\tau_1$/5 | 4$\tau_1$/5 | $\tau_1$ |
|---|---|---|---|---|---|---|
| -1.000E+00 | 1.06369836E+00 | 8.58242289E-01 | 6.35970510E-01 | 4.13094529E-01 | 1.90467611E-01 | 0.00000000E+00 |
| -9.000E-01 | 9.53090077E-01 | 8.69389789E-01 | 6.47424656E-01 | 4.24549088E-01 | 2.01733846E-01 | 0.00000000E+00 |
| -8.000E-01 | 9.54076470E-01 | 8.80528187E-01 | 6.58878508E-01 | 4.36004258E-01 | 2.13094528E-01 | 0.00000000E+00 |
| -7.000E-01 | 8.82541842E-01 | 8.91568315E-01 | 6.70332018E-01 | 4.47459699E-01 | 2.24510135E-01 | 0.00000000E+00 |
| -6.000E-01 | 8.24712324E-01 | 9.02556262E-01 | 6.81785223E-01 | 4.58915261E-01 | 2.35953990E-01 | 0.00000000E+00 |
| -5.000E-01 | 7.72605679E-01 | 9.13497493E-01 | 6.93238122E-01 | 4.70370875E-01 | 2.47410721E-01 | 0.00000000E+00 |
| -4.000E-01 | 7.11438505E-01 | 9.24375836E-01 | 7.04690673E-01 | 4.81826511E-01 | 2.58872714E-01 | 0.00000000E+00 |
| -3.000E-01 | 6.40310565E-01 | 9.35178929E-01 | 7.16142813E-01 | 4.93282152E-01 | 2.70336510E-01 | 0.00000000E+00 |
| -2.000E-01 | 5.58481727E-01 | 9.45894466E-01 | 7.27594473E-01 | 5.04737793E-01 | 2.81800546E-01 | 0.00000000E+00 |
| -1.000E-01 | 4.58734040E-01 | 9.56509437E-01 | 7.39045581E-01 | 5.16193429E-01 | 2.93264125E-01 | 0.00000000E+00 |
| 0.000E+00 | 2.51582455E-01 | 9.67014307E-01 | 7.50496082E-01 | 5.27649059E-01 | 3.04726966E-01 | 3.92638589E-02 |
| 1.000E-01 | 0.00000000E+00 | 9.77407148E-01 | 7.61945957E-01 | 5.39104682E-01 | 3.16188998E-01 | 7.20390691E-02 |
| 2.000E-01 | 0.00000000E+00 | 9.87701480E-01 | 7.73395273E-01 | 5.50560298E-01 | 3.27650253E-01 | 8.89489741E-02 |
| 3.000E-01 | 0.00000000E+00 | 9.97940900E-01 | 7.84844259E-01 | 5.62015910E-01 | 3.39110806E-01 | 1.04140400E-01 |
| 4.000E-01 | 0.00000000E+00 | 1.00822662E+00 | 7.96293473E-01 | 5.73471520E-01 | 3.50570747E-01 | 1.18383923E-01 |
| 5.000E-01 | 0.00000000E+00 | 1.01877096E+00 | 8.07744108E-01 | 5.84927135E-01 | 3.62030161E-01 | 1.31990433E-01 |
| 6.000E-01 | 0.00000000E+00 | 1.03000170E+00 | 8.19198648E-01 | 5.96382767E-01 | 3.73489122E-01 | 1.45139314E-01 |
| 7.000E-01 | 0.00000000E+00 | 1.04278189E+00 | 8.30662263E-01 | 6.07838443E-01 | 3.84947681E-01 | 1.57957743E-01 |
| 8.000E-01 | 0.00000000E+00 | 1.05891057E+00 | 8.42146047E-01 | 6.19294219E-01 | 3.96405868E-01 | 1.70543344E-01 |
| 9.000E-01 | 0.00000000E+00 | 1.08265042E+00 | 8.53675508E-01 | 6.30750225E-01 | 4.07863688E-01 | 1.82954318E-01 |
| 1.000E+00 | 0.00000000E+00 | 1.26059598E+00 | 8.65364223E-01 | 6.42206808E-01 | 4.19321107E-01 | 1.95231197E-01 |

Finally, Table IIb gives the 8-place benchmark for uniform spatial edits and 21 directional edits, we see two discrepancies in the exiting intensities.

Therefore, from this rather incomplete survey of benchmark comparisons, we are at least confident in 7 places for 1$^{st}$/DOM/RM.

## CONCLUSIONS

A new response matrix algebraic solution, based on the first order solution of the 1D linear transport equation for general anisotropic scattering has been presented. We compare the matrix exponential form, the most fundamental of all solutions of the transport equation, to the second order solution of RM/DOM for the 300-term CloudC1 Legendre polynomial scattering kernel. The two solutions differ in the choice of the linearly independent solutions to represent the solution to the homogeneous equation. The two different forms agreed to 8 significant digits. A primary reason for the high precision of the scaled solution comes from Wynn-epsilon acceleration [5].

## REFERENCES

[1]. Ganapol, B.D., *The response matrix discrete ordinates solution to the 1D radiative transfer equation*. Journal of Quantitative Spectroscopy & Radiative Transfer, **154**, pp. 72–90 (2015).
[2]. Moler ,C. and Van Loan, C*., Nineteen Dubious Ways to Compute the Exponential of a Matrix, Twenty-Five Years Later,* SIAM Review, **V45**,#1, pp. 3 -47, 2003.
[3] Siewert, C.E., *A concise and accurate solution to Chandrasekhar's basic problem in radiative transfer*, Journal of Quantitative Spectroscopy & Radiative Transfer, **64** pp. 109-130, (2000).
[4] Garcia, R.D.M. and Siewert, C.E., *Benchmark results in radiative transfer*, TTSP, **14**, 437-483, (1985).
[5] Wynn, P., *On a device for computing the em(Sn) transformation*, Math., Tables Aids Comput. **10**, pp. 91–96, (1956).